\documentclass[a4paper,10pt]{article}

\usepackage[latin1]{inputenc}
\usepackage{amsfonts}
\usepackage{amsmath}
\usepackage{amssymb}
\usepackage{amsthm}
\usepackage[T1]{fontenc}
\usepackage[dvips]{graphicx}

\author{Robert J Taggart}
\title{Pointwise Convergence for Semigroups in Vector-valued $L^p$ Spaces}
\date{18 February 2008}

\newtheorem{theorem}{Theorem}[section]
\newtheorem{corollary}[theorem]{Corollary}
\newtheorem{lemma}[theorem]{Lemma}

\theoremstyle{definition}
\newtheorem{definition}[theorem]{Definition}

\theoremstyle{remark}

\newcommand{\CC}{\mathbb{C}}
\newcommand{\R}{\mathbb{R}}
\newcommand{\Z}{\mathbb{Z}}
\newcommand{\N}{\mathbb{N}}
\newcommand{\Q}{\mathbb{Q}}

\newcommand{\bs}{\backslash}
\newcommand{\B}{\mathcal{B}}

\newcommand{\Hi}{\mathcal{H}}
\newcommand{\U}{\mathcal{U}}
\newcommand{\real}{\mathrm{Re}}
\newcommand{\im}{\mathrm{Im}}
\newcommand{\norm}[1]{\left\Vert{#1}\right\Vert}
\newcommand{\wt}[1]{\widetilde{#1}}
\newcommand{\dd}{\mathrm{d}}

\newcommand{\bT}{\mathbf{T}}

\def\ip<#1,#2>{\left\langle#1,\,#2\right\rangle}

\begin{document}

\maketitle
\begin{abstract}
\noindent Suppose that $\{T_t:t\geq0\}$ is a symmetric diffusion semigroup on $L^2(X)$ and denote by $\{\wt{T}_t:t\geq0\}$ its tensor product extension to the Bochner space $L^p(X,\B)$, where $\B$ belongs to a certain broad class of UMD spaces. We prove a vector-valued version of the Hopf--Dunford--Schwartz ergodic theorem and show that this extends to a maximal theorem for analytic continuations of $\{\wt{T}_t:t\geq0\}$ on $L^p(X,\B)$. As an application, we show that such continuations exhibit pointwise convergence.
\end{abstract}
\section{Introduction}

The goal of this paper is to show that two classical results about symmetric diffusion semigroups, which go back to E. M. Stein's monograph \cite{eS70}, can be extended to the setting of vector-valued $L^p$ spaces.

Suppose throughout that $(X,\mu)$ is a positive $\sigma$-finite measure space.

\begin{definition} Suppose that $\{T_t:t\geq0\}$ is a semigroup of operators on $L^2(X)$. We say that
\begin{enumerate}
\item[(a)] the semigroup $\{T_t:t\geq0\}$ satisfies the \textit{contraction property} if
\begin{equation}\label{eq:contraction}
\norm{T_tf}_q\leq\norm{f}_q\qquad\forall f\in L^2(X)\cap L^q(X)
\end{equation}
whenever $t\geq0$ and $q\in[1,\infty]$; and
\item[(b)] the semigroup $\{T_t:t\geq0\}$ is a \textit{symmetric diffusion semigroup} if it satisfies the contraction property and if $T_t$ is selfadjoint on $L^2(X)$ whenever $t\geq0$.
\end{enumerate}
\end{definition}

It is well known that if $1\leq p<\infty$ and $1\leq q\leq\infty$ then $L^q(X)\cap L^p(X)$ is dense in $L^p(X)$. Hence, if a semigroup $\{T_t:t\geq0\}$ acting on $L^2(X)$ has the contraction property then each $T_t$ extends uniquely to a contraction of $L^p(X)$ whenever $p\in[1,\infty)$. By abuse of notation, we shall also denote by $\{T_t:t\geq0\}$ the unique semigroup extension which acts on $L^p(X)$.

The class of symmetric diffusion semigroups is widely used in applications and includes the Gaussian and Poisson semigroups on $L^2(\R^n)$. Despite the simplicity of the axioms defining this class, symmetric diffusion semigroups have a rich theory. 
For example, if $\{T_t:t\geq0\}$ is a symmetric diffusion semigroup on $L^2(X)$ then the semigroup can also be continued analytically to sectors of the complex plane. To be precise, given a positive angle $\psi$, let $\Gamma_{\psi}$ denote the cone $\{z\in\CC:|\arg z|<\psi\}$ and $\overline{\Gamma}_{\psi}$ its closure. We shall denote the interval $[0,\infty)$ by $\overline{\Gamma}_0$. Using spectral theory and complex interpolation, Stein proved the following result.

\begin{theorem}[Stein \cite{eS70}]\label{th:analytic continuation} Suppose that $1<p<\infty$,
\[\psi/\pi=1/2-|1/p-1/2|>0,\]
and $\{T_t:t\in\R\}$ is a symmetric diffusion semigroup on $L^2(X)$. Then $\{T_t:t\geq0\}$ extends uniquely to a semigroup $\{T_z:z\in\overline{\Gamma}_{\psi}\}$ of contractions on $L^p(X)$ such that the operator-valued function $z\mapsto T_z$ is holomorphic in $\Gamma_{\psi}$ and weak operator topology continuous in $\overline{\Gamma}_{\psi}$.
\end{theorem}

We now recall two results of M. Cowling \cite{mC83}, developing the fundamental work of Stein \cite{eS70}, which the current paper generalises. The first is a useful technical tool. For $f$ in $L^p(X)$, define the maximal function $\mathcal{M}^{\psi}f$ by
\[\mathcal{M}^{\psi}f=\sup\{|T_zf|:z\in\overline{\Gamma}_{\psi}\}.\]
The maximal theorem, stated below, says that the maximal function operator $\mathcal{M}^{\psi}$ is bounded on $L^p(X)$.

\begin{theorem}[Stein--Cowling \cite{mC83}]\label{th:scalar maximal}
Suppose that $1<p<\infty$ and that
\[0\leq \psi/\pi<1/2-|1/p-1/2|.\]
If $\{T_z:z\in\overline{\Gamma}_{\psi}\}$ is the semigroup on $L^p(X)$ given by Theorem \ref{th:analytic continuation} then there is a positive constant $C$ such that
\[\norm{\mathcal{M}^{\psi}f}_p\leq C\norm{f}_p\qquad\forall f\in L^p(X).\]
\end{theorem}

The maximal theorem allows one to deduce a pointwise convergence result for the semigroup $\{T_z:z\in\overline{\Gamma}_{\psi}\}$.

\begin{corollary}[Stein--Cowling \cite{mC83}]\label{cor:sv ptwise}  Assume the hypotheses of Theorem \ref{th:scalar maximal}. If $f\in L^p(X)$ then $(T_zf)(x)\to f(x)$ for almost every $x$ in $X$ as $z$ tends to $0$ in $\overline{\Gamma}_{\psi}$.
\end{corollary}

The earliest form of the maximal theorem appeared in Stein \cite[p.~73]{eS70} for the case when $\psi=0$. From this Stein deduced the pointwise convergence of $T_tf$ to $f$ as $t\to 0^+$. Using a simpler approach, Cowling \cite{mC83} extended Stein's result to semigroups $\{T_z:z\in\overline{\Gamma}_{\psi}\}$, holomorphic in the sector $\Gamma_{\psi}$, without additional hypotheses. Given $z\in\overline{\Gamma}_{\psi}$, Cowling's strategy was to decompose the operator $T_z$ into two parts:
\begin{equation}\label{eq:decomposition}
T_zf=\frac{1}{t}\int_0^te^{-sL}f\,\dd s+\Big[e^{-zL}f-\frac{1}{t}\int_0^te^{-sL}f\,\dd s\Big],
\end{equation}
where $t=|z|$ and $-L$ is the generator of the semigroup. The $L^p$ norm of the first term on the right-hand side can be controlled by the Hopf--Dunford--Schwartz ergodic theorem. A clever use of the Mellin transform allows the bracketed terms to be controlled by bounds on the imaginary powers of $L$.

The chief contribution of this paper is to observe that, under certain assumptions, this argument may be adapted to the setting of $L^p$ spaces of Banach-space-valued functions. Several other results contained in Stein's monograph \cite{eS70} have already been pushed in this direction (see, for example, \cite{qX98}, \cite{MTX06} and \cite{tH07}). In a broader context, there has been much recent interest in operators which act on such spaces, particularly when the Banach space has the so-called UMD property (see the ground breaking work of J. Bourgain \cite{jB83} and D. Burkholder \cite{dB83}). Developments which are perhaps most pertinent to our results include studies on bounded imaginary powers of operators (of which the article \cite{DV87} of G. Dore and A. Venni is now a classic), $H^{\infty}$-functional calculi for sectorial operators (see especially the paper \cite{CDMY96} of A. McIntosh and his collaborators) and maximal $L^p$-regularity (see L. Weis \cite{lW01} and the references therein). The article \cite{KW04} of P. Kunstmann and L. Weis gives an excellent exposition of the interplay between these motifs in the vector-valued setting as well as an extensive bibliography detailing the key contributions made to the field over the last two decades.

Suppose that $\B$ is a (complex) Banach space and let $L^p(X,\B)$ denote the Bochner space of $\B$-valued $p$-integrable functions on $X$. Given a symmetric diffusion semigroup $\{T_t:t\geq0\}$ on $L^2(X)$, its tensor product extension $\{\wt{T}_t:t\geq0\}$ to $L^p(X,\B)$ exists by the contraction property (see Section \ref{s:vv_ext}). If $\{\wt{T}_t:t\geq0\}$ can be continued analytically to some sector $\Gamma_{\psi+\epsilon}$, where $0<\psi<\pi/2$ and $\epsilon$ is a (sufficiently) small positive number, then denote this continuation by $\{\wt{T}_z:z\in\Gamma_{\psi+\epsilon}\}$. If such a continuation does not exist, we take $\psi$ to be $0$.  Given any function $F$ in $L^p(X,\B)$, one defines the maximal function $\mathcal{M}^{\psi}_{\B}F$ by
\begin{equation}\label{eq:def of maximal function}
\mathcal{M}^{\psi}_{\B}F=\sup\{|\wt{T}_zF|_{\B}:z\in\overline{\Gamma}_{\psi}\}.
\end{equation}
The theorem below is the main result of this paper.

\begin{theorem}\label{th:main semigroup theorem}
Suppose that $(X,\mu)$ is a $\sigma$-finite measure space and that $\{T_t:t\geq0\}$ is a symmetric diffusion semigroup on $L^2(X)$. Suppose also that $\B$ is a Banach space isomorphic to a closed subquotient of a complex interpolation space $(\Hi,\U)_{[\theta]}$, where $\Hi$ is a Hilbert space, $\U$ is a UMD space and $0<\theta<1$. If $1<p<\infty$, $|2/p-1|<\theta$ and
\[0\leq\psi<\frac{\pi}{2}(1-\theta)\]
then 
\begin{enumerate}
\item[(a)] $\{\wt{T}_t:t\geq0\}$ has a bounded analytic continuation to the sector $\Gamma_{\psi}$ in $L^p(X,\B)$,
\item[(b)] there is a positive constant $C$ such that
\[\norm{\mathcal{M}^{\psi}_{\B}F}_{L^p(X)}\leq C\norm{F}_{L^p(X,\B)}\qquad\forall F\in L^p(X,\B),\]
and
\item[(c)] if $F\in L^p(X,\B)$ then $\wt{T}_zF(x)$ converges to $F(x)$ for almost every $x$ in $X$ as $z$ tends to $0$ in the sector $\overline{\Gamma}_{\psi}$.
\end{enumerate}
\end{theorem}

It is noteworthy that the class of Banach spaces $\B$ satisfying the interpolation hypothesis of Theorem \ref{th:main semigroup theorem} is a subset of those Banach spaces possessing the UMD property. It includes those classical Lebesgue spaces, Sobolev spaces and Schatten--von Neumann ideals that are reflexive. The reader is directed to Section \ref{s:BIP} for further remarks on these spaces.

The structure and content of the rest of this paper is as follows. Section \ref{s:vv_ext} presents some standard results on tensor product extensions of operators to vector-valued $L^p$ spaces. For example, it is well-known that such extensions exist for semigroups with the contraction property and consequently these semigroups are \textit{subpositive}. In Section \ref{s:subpositivity} we prove a stronger result; namely that, whenever $1\leq p<\infty$, every measurable semigroup $\{T_t:t\geq0\}$ on $L^2(X)$ satisfying the contraction property is dominated on $L^p(X)$ by a measurable positive semigroup with the contraction property. This result allows us to easily deduce, in Section \ref{s:HDS}, a vector-valued version of the Hopf--Dunford--Schwartz ergodic theorem.
 
Parts (a) and (b) of Theorem \ref{th:main semigroup theorem} are proved in Sections \ref{s:maximal theorem} and \ref{s:BIP}. Following techniques used in \cite{mC83}, we begin by proving a maximal theorem for the tensor prouct extension $\{\wt{T}_t:t\geq0\}$ to $L^p(X,\B)$ of a strongly continuous semigroup $\{T_t:t\geq0\}$ satisfying the contraction property. Here we assume that $1<p<\infty$ and $\B$ is any Banach space, provided that the generator $-\wt{L}$ of the $\B$-valued extension has bounded imaginary powers on $L^p(X,\B)$ with a power angle less than $\pi/2-\psi$. Section \ref{s:BIP} discusses circumstances under which this condition holds. In general, it is necessary that $\B$ has the UMD property. Moreover, by exploiting the subpositivity of $\{T_t:t\geq0\}$ and adapting arguments of M. Hieber and J. Pr\"uss \cite{HP98}, we show that if $\B$ has the UMD property then $\wt{L}$ has an $H^{\infty}$-functional calculus. This, along with spectral theory (where the self-adjointness of each operator $T_t$ is imposed) and interpolation, allows us to remove the bounded imaginary powers hypothesis at the cost of restricting the class of Banach spaces $\B$ for which the maximal theorem is valid.

In Section \ref{s:proof of main result} we show that that the pointwise convergence of $\{\wt{T}_z:z\in\overline{\Gamma}_{\psi}\}$ is easily deduced from the pointwise convergence of $\{T_z:z\in\overline{\Gamma}_{\psi}\}$ and the maximal theorem. This completes the proof of Theorem \ref{th:main semigroup theorem}.

\section{Vector-valued extensions of contraction semigroups}\label{s:vv_ext}

Suppose that $\B$ is a (complex) Banach space with norm $|\cdot|_{\B}$ and that $(X,\mu)$ is a $\sigma$-finite measure space. We also assume throughout this section that $p\in[1,\infty)$. Denote by $L^p(X,\B)$ the Bochner space of all $\B$-valued measurable functions $F$ on $X$ satisfying
\[\norm{F}_p:=\Big(\int_X|F(x)|_{\B}^p\,\dd \mu(x)\Big)^{1/p}<\infty.\]
(As is customary, we will not distinguish between equivalence classes of functions and members of each equivalence class.)
Let $L^p(X)\otimes\B$ denote the set of all finite linear combinations of $\B$-valued functions of the form $uf$, where $u\in\B$ and $f\in L^p(X)$.
It is known that this set is dense in $L^p(X,\B)$. Many operators acting on scalar-valued function spaces can be extended to act on $\B$-valued function spaces in the following canonical way.

\begin{definition}
Suppose that $T$ is a bounded operator on $L^p(X)$. If $I_{\B}$ denotes the identity operator on $\B$ then define the tensor product $T\otimes I_{\B}$ on $L^p(X)\otimes\B$ by
\[T\otimes I_{\B}\left(\sum_{k=1}^nu_kf_k\right)=\sum_{k=1}^nu_kTf_k\]
whenever $n\in\Z^+$, $u_k\in\B$, $f_k\in L^p(X)$ and $k=1,\ldots,n$.
We say that a bounded operator $\wt{T}:L^p(X,\B)\to L^p(X,\B)$ is a \textit{$\B$-valued extension of $T$} if $\wt{T}=T\otimes I_{\B}$ on $L^p(X)\otimes\B$. In this case $\wt{T}$ is also called a \textit{tensor product extension of $T$ to $L^p(X,\B)$}.
\end{definition}

If it exists, a $\B$-valued extension $\wt{T}$ of $T$ is necessarily unique by the density of $L^p(X)\otimes\B$ in $L^p(X,\B)$. 

It is well-known that if an operator $T$ on $L^2(X)$ extends to a contraction on $L^q(X)$ for all $q$ in $[1,\infty]$ then $T\otimes I_{\B}$ extends to a contraction $\wt{T}$ on $L^p(X,\B)$. (This is not hard to show if $p=1$; for other values of $p$ the result can be deduced by duality and interpolation.) Consequently, any semigroup $\{T_t:t\geq0\}$ on $L^2(X)$ with the contraction property extends to a semigroup $\{\wt{T}_t:t\geq0\}$ of contractions on $L^p(X,\B)$. Moreover, if $\{T_t:t\geq0\}$ is strongly continuous on $L^p(X)$ then $\{\wt{T}_t:t\geq0\}$ is strongly continuous on $L^p(X,\B)$. This is a consequence of the following lemma.

\begin{lemma}\label{lem:strong continuity}
Suppose that $\{S_t:t\geq0\}$ is a family of bounded operators on $L^p(X)$ with a $\B$-valued extension $\{\wt{S}_t:t\geq0\}$ to $L^p(X,\B)$. If the mapping $t\mapsto S_t$ is strongly continuous and $\{\wt{S}_t:t\geq0\}$ is locally (with respect to $t$) uniformly bounded in norm then the mapping $t\mapsto\wt{S}_t$ is strongly continuous.
\end{lemma}

\begin{proof}
Suppose that $F\in L^p(X,\B)$. Then we may approximate $F$ by a function $G$ of the form $\sum_{k=1}^nu_kg_k$, where $n\in\N$, $u_k\in\B$ and $g_k\in L^p(X)$. Since the map $t\mapsto S_t$ is strongly continuous, the map $t\mapsto \wt{S}_tG$ is continuous, and by a standard $3\epsilon$ argument, so is the map $t\mapsto \wt{S}_tF$. 
\end{proof}

Suppose that $\{T_t:t\geq0\}$ is a strongly continuous semigroup on $L^2(X)$ with the contraction property. The generator $B$ of tensor extension $\{\wt{T}_t:t\geq0\}$ to $L^p(X,\B)$ is given, as usual, by
\[BF=\lim_{t\to0^+}\frac{\wt{T}_tF-F}{t}\]
for all $F$ in $L^p(X,\B)$ for which the limit exists. The collection of such $F$ is called the domain of $B$. Let $-L$ denote the generator of $\{T_t:t\geq0\}$. It is easy to show that $\mathrm{Dom}(L)\otimes\B\subseteq\mathrm{Dom}(B)$ and that $B=-L\otimes I_{\B}$ on $\mathrm{Dom}(L)\otimes\B$. Therefore we shall denote $B$ by $-\wt{L}$.

Bounded operators with vector-valued extensions can be characterised in terms of \textit{subpositivity}, a property that proves useful in subsequent sections.

\begin{definition}\label{def:subpositive} Suppose that $T$ is a linear operator on $L^p(X)$. We say that
\begin{enumerate}
\item[(a)] $T$ is \textit{positive} if $Tf\geq0$ whenever $f\geq0$ for $f$ in $L^p(X)$;
\item[(b)] $T$ is \textit{subpositive} if there exists a bounded positive operator $S$ on $L^p(X)$ such that $|Tf|\leq S|f|$ whenever $f\in L^p(X)$, in which case we also say that $T$ is \textit{dominated} by $S$.
\end{enumerate}
\end{definition}

If $R$ is an operator on $L^p(X)$ then define $\overline{R}$ by the formula $\overline{R}f=\overline{R\bar{f}}$ whenever $f\in L^p(X)$, and define $\real(R)$ by $(R+\overline{R})/2$.

\begin{lemma}\label{lem:subpositivity characterisations}
Suppose that $T$ is a bounded operator on $L^p(X)$. Then the following are equivalent.
\begin{enumerate}
\item[(a)] The operator $T$ is a subpositive contraction on $L^p(X)$.
\item[(b)] For any finite subset $\{f_k\}_{k=1}^n$ of $L^p(X)$,
\[\norm{\sup_k|Tf_k|}_{L^p(X)}\leq\norm{\sup_k|f_k|\,}_{L^p(X)}.\]
\item[(c)] The tensor product $T\otimes I_{\B}$ extends to a contraction $\wt{T}$ on $L^p(X,\B)$.
\item[(d)] There exists a positive contraction $S$ such that $S+\real(e^{i\theta}T)$ is positive whenever $\theta\in\R$.
\end{enumerate}
\end{lemma}

The equivalence of (a), (b) and (c) are well-known (see, for example, \cite{gP94} and \cite{vP83}) and will be used in Section \ref{s:maximal theorem}. Statement (d) is the definition of subpositive contractivity given by R. Coifman, R. Rochberg and G. Weiss \cite[p.~54]{CRW77}. The fact that (a) implies (d) is easy to establish and is used to prove Theorem \ref{th:H infty functional calculus}. The converse is harder to prove (see, for example, \cite[Section 2.1]{rT08}) but will not be needed for the results of this paper.

\section{Subpositivity for contraction semigroups}\label{s:subpositivity}

One of the consequences of Lemma \ref{lem:subpositivity characterisations} is that any semigroup on $L^2(X)$ that enjoys the contraction property extends to a semigroup of subpositive contractions on $L^p(X)$ whenever $1\leq p<\infty$. The goal of this section is to prove a stronger result needed in Section \ref{s:HDS}, namely that every semigroup on $L^2(X)$ with the contraction property is, when extended to a semigroup on $L^p(X)$ for $p$ in $[1,\infty)$, dominated by a positive contraction semigroup on $L^p(X)$.

We begin with a few preliminaries. Suppose that $1\leq p<\infty$ and $T$ is a bounded linear operator on $L^p(X)$. If $1\leq q<\infty$ and $\norm{Tf}_q\leq C\norm{f}_q$ for all $f$ in $L^q(X)\cap L^p(X)$ then $T$ has a unique bounded linear extension acting on $L^q(X)$. By abuse of notation we will also denote this extension by $T$.

We say that a family of operators $\{T_t:t\geq0\}$ is \textit{(strongly) measurable} on $L^p(X)$ if, for every $f$ in $L^p(X)$, the $L^p(X)$-valued map $t\mapsto T_tf$ is measurable with respect to Lebesgue measure on $[0,\infty)$. The family is said to be \textit{weakly measurable} if the complex-valued map $t\mapsto\ip<T_tf,g>$ is measurable with respect to Lebesgue measure on $[0,\infty)$ whenever $f\in L^p(X)$ and $g\in L^{p'}(X)$. Here, $p'$ denotes the conjugate exponent of $p$, given by $1/p+1/p'=1$. If $1\leq p<\infty$ then $L^p(X)$ is a separable Banach space and hence strong measurability and weak measurability coincide by the Pettis measurability theorem (see \cite[Theorem III.6.11]{DS58}).

We now state the main result of this section.

\begin{theorem}\label{th:contraction implies subpositive}
Suppose that $\{T_t:t\geq0\}$ is a semigroup on $L^2(X)$ satisfying the contraction property. Then there exists a positive semigroup $\{S_t:t\geq0\}$ on $L^2(X)$, satisfying the contraction property, such that
\[|T_tf|\leq S_t|f|\qquad\forall f\in L^p(X)\]
whenever $1\leq p<\infty$ and $t\geq0$.
If $\{T_t:t\geq0\}$ is a measurable semigroup on $L^2(X)$ then $\{S_t:t\geq0\}$ extends to a measurable semigroup on $L^p(X)$ whenever $1\leq p<\infty$.
\end{theorem}

We shall prove the theorem via a sequence of lemmata. The final stage of the proof draws heavily on the work of Y. Kubokawa \cite{yK75} and C. Kipnis \cite{cK74}, who independently proved a similar result for $L^1$ contraction semigroups.

Denote by the $L^p_+(X)$ the set of nonnegative functions in $L^p(X)$.

\begin{lemma}\label{lem:subpositive on all Lp}
Suppose that $1<p<\infty$. Assume also that $T$ and $S$ are bounded operators on $L^1(X)$ such that $\norm{Tf}_p\leq\norm{f}_p$ and $\norm{Sf}_p\leq\norm{f}_p$ whenever $f\in L^1(X)\cap L^p(X)$. If $S$ is positive and dominates $T$ on $L^1(X)$ then
\[|Tf|\leq S|f|\qquad\forall f\in L^p(X).\]
\end{lemma}

\begin{proof}
Assume the hypotheses of the lemma and suppose that $f\in\ L^p(X)$. Since $L^1(X)\cap L^p(X)$ is dense in $L^p(X)$ there is a sequence $\{f_n\}_{n\in\N}$ in $L^1(X)\cap L^p(X)$ such that $f_n\to f$ in $L^p(X)$. By continuity, $|Tf_n|\to |Tf|$ in $L^p(X)$ and similarly $S|f_n|\to S|f|$ in $L^p(X)$. Moreover, $|Tf_n|\leq S|f_n|$ for all $n$. If $g_n=S|f_n|-|Tf_n|$ and $g=S|f|-|Tf|$ then each $g_n$ is nonnegative and $g_n\to g$ in $L^p(X)$. Now $L^p_+(X)$ is a closed subset of $L^p(X)$, so $g\geq0$ and the proof is complete.
\end{proof}

\begin{lemma}\label{lem:measurable for all p}
Suppose that $p$ and $q$ both lie in the interval $[1,\infty)$ and that $\{T_t:t\geq0\}$ is a family of operators on $L^2(X)$ satisfying the contraction property. If $\{T_t:t\geq0\}$ is measurable on $L^p(X)$ then $\{T_t:t\geq0\}$ is measurable on $L^q(X)$.
\end{lemma}

\begin{proof}
Assume the hypotheses and suppose that $f\in L^q(X)$ and $g\in L^{q'}(X)$. It suffices to show that the map $\phi:[0,\infty)\to\CC$, defined by
\[\phi(t)=\ip<T_tf,g>,\]
is measurable. Choose any sequence $\{X_n\}_{n\in\N}$ of measurable sets satisfying $X_n\subset X_{n+1}$ and $\cup_{n\in\N}X_n=X$. By carefully choosing a sequence $\{g_n\}_{n\in\N}$ contained in $L^{q'}(X)\cap L^{p'}(X)$ which converges in $L^{q'}(X)$ to $g$, and a sequence $\{f_n\}_{n\in\N}$ contained in $L^{q}(X)\cap L^{p}(X)$ which converges in $L^q(X)$ to  $f$, one can show that the sequence $\{\phi_n\}_{n\in\N}$, given by
\[\phi_n(t)=\ip<T_tf_n,1_{X_n}g_n>,\]
converges pointwise to $\phi$. Since each $\phi_n$ is measurable, $\phi$ is also.
\end{proof}

\begin{lemma}\label{lem:linear modulus}
Suppose that $T$ is a linear contraction on $L^1(X)$ with the property that $\norm{Tf}_q\leq\norm{f}_q$ whenever $f\in L^q(X)\cap L^1(X)$ and $1\leq q\leq \infty$. Then there is a unique bounded linear positive operator $\bT$ on $L^1(X)$ such that
\begin{enumerate}
\item[(a)] the operator norms of $T$ and $\bT$ on $L^1(X)$ are equal,
\item[(b)] $\norm{\bT f}_q\leq\norm{f}_q$ whenever $f\in L^q(X)\cap L^1(X)$ and $1\leq q\leq \infty$,
\item[(c)] $|Tf|\leq\bT|f|$ whenever $f\in L^p(X)$ and $1\leq p\leq \infty$, and
\item[(d)] $\bT f=\sup\{|Tg|:g\in L^1(X), |g|\leq f\}$ whenever $f\in L^1_+(X)$.
\end{enumerate}
\end{lemma}

\begin{proof}
For the existence of a unique operator $\bT$ satisfying properties (a), (c) (for the case when $p=1$) and (d), see, for example, \cite[Theorem 4.1.1]{uK85}. Property (b) holds by \cite[Lemma VIII.6.4]{DS58}. We can now deduce property (c), for the case when $1<p<\infty$, from Lemma \ref{lem:subpositive on all Lp}.
\end{proof}

The operator $\bT$ introduced in the lemma is called the \textit{linear modulus} of $T$. If $\{T_t:t\geq0\}$ is a bounded semigroup on $L^1(X)$ then $\bT_{s+t}\leq\bT_s\bT_t$ for all nonnegative $s$ and $t$. However, equality may not hold and thus the family $\{\bT_t:t\geq0\}$ of bounded positive operators will not, in general, be a semigroup. Nevertheless, Kubokawa \cite{yK75} and Kipnis \cite{cK74} (see \cite[Theorems 4.1.1 and 7.2.7]{uK85} for a more recent exposition) showed that the linear modulus $\bT_t$ could be used to construct a positive semigroup $\{S_t:t\geq0\}$, known as the \textit{modulus semigroup}, which dominates $\{T_t:t\geq0\}$. The following proof uses this construction.

\medskip

\noindent{\itshape\rmfamily Proof of Theorem \ref{th:contraction implies subpositive}.}
Assume the hypothesis of Theorem \ref{th:contraction implies subpositive} and suppose that $t>0$. Let $\mathcal{D}_t$ denote the family of all finite subdivisions $(s_i)$ of $[0,t]$ satisfying
\[0=s_0<s_1<s_2<\ldots<s_n=t.\]
If $\mathbf{s}=(s_i)$ and $\mathbf{s}'=(s_j')$ are two elements of $\mathcal{D}_t$ then we write $\mathbf{s}<\mathbf{s}'$ whenever $\mathbf{s}'$ is a refinement of $\mathbf{s}$. With this partial order $\mathcal{D}_t$ is an increasingly filtered set. For $f$ in $L^1_+(X)$, put
\[\Phi(\mathbf{s},f)=\bT_{s_1}\bT_{s_2-s_1}\ldots\bT_{s_n-s_{n-1}}f,\]
where $\bT_{\alpha}$ is the linear modulus of $T_{\alpha}$ whenever $\alpha\geq0$.
It follows from $\bT_{\alpha+\beta}\leq\bT_{\alpha}\bT_{\beta}$ that $\mathbf{s}<\mathbf{s}'$ implies $\Phi(\mathbf{s},f)\leq\Phi(\mathbf{s}',f)$. Since the operator $\bT_{\alpha}$ is contraction whenever $\alpha\geq0$, we have $\norm{\Phi(\mathbf{s},f)}_1\leq\norm{f}_1$. We now define $S_t$ on $L^1_+(X)$ by
\[S_tf=\sup\{\Phi(\mathbf{s},f):\mathbf{s}\in\mathcal{D}_t\}.\]
Note that
\[\sup\{\Phi(\mathbf{s},f):\mathbf{s}\in\mathcal{D}_t\}=\lim_{\mathbf{s}\in\mathcal{D}_t}\Phi(\mathbf{s},f)\]
so $S_t$ is well-defined by the monotone convergence theorem for increasingly filtered families.

It is easy to check that $S_t(f+g)=S_tf+S_tg$ and $S_t(\lambda f)=\lambda S_tf$ whenever $f$ and $g$ belong to $L^1_+(X)$ and $\lambda\geq0$. Moreover, $\norm{S_tf}_1\leq\norm{f}_1$ if $f\in L^1_+(X)$. Therefore $S_t$ can now be defined for all $f$ in $L^1(X)$ as a linear contraction of $L^1(X)$. We define $S_0$ as the identity operator on $L^1(X)$.

We now show that $\{S_t:t\geq0\}$ is a semigroup. Suppose that $t$ and $t'$ are both positive. If
\[0=s_0<s_1<s_2<\ldots<s_n=t\]
and
\[0=s_0'<s_1'<s_2'<\ldots<s_n'=t'\]
form subdivisions of $[0,t]$ and $[0,t']$ then
\[0=s_0<s_1<s_2<\ldots<s_n=s_n+s_0'<s_n+s_1'<s_n+s_2'<\ldots<s_n+s_n'\]
forms a subdivision of $[0,t+t']$. Conversely every subdivision of $[0,t+t']$ which is fine enough to contain $t$ is of this form. This yields $S_{t+t'}=S_tS_{t'}$.

By Lemma \ref{lem:linear modulus} (b), it is easy to check that $\{S_t:t\geq0\}$ extends to a contraction semigroup on $L^2(X)$ which satisfies the contraction property. Moreover, the construction shows that $|T_tf|\leq S_t|f|$ whenever $f\in L^1(X)$ and $t\geq0$. By an application of Lemma \ref{lem:subpositive on all Lp}, we deduce that each $T_t$ is dominated by $S_t$ on $L^p(X)$ whenever $1\leq p<\infty$.

It remains to show that if $\{T_t:t\geq0\}$ is measurable on $L^2(X)$ then $\{S_t:t\geq0\}$ is measurable on $L^p(X)$ whenever $1\leq p<\infty$. In view of Lemma \ref{lem:measurable for all p}, we can assume that $\{T_t:t\geq0\}$ is measurable on $L^1(X)$ and it suffices to show that $\{S_t:t\geq0\}$ is measurable on $L^1(X)$. Fix $f$ in $L^1(X)$ and define $\phi:[0,\infty)\to L^1(X)$ by $\phi(t)=S_tf$. We will construct a sequence $\{\phi_n\}_{n\in\N}$ of measurable functions converging pointwise to $\phi$, completing the proof.

Since $f$ can be decomposed as a linear combination of four nonnegative functions (the positive and negative parts of $\real(f)$ and $\im(f)$) and each $S_t$ is linear, we can assume, without loss of generality, that $f\geq0$.

When $t>0$, $n\in\N$ and $m$ is the smallest integer such that $m\geq 2^nt$, let $\mathbf{s}(n,t)$ denote the subdivision $(s_k(n,t))_{k=0}^m$ of $[0,t]$ given by
\[s_k(n,t)=\begin{cases}
k2^{-n} & \mbox{if }k=0,1,\ldots,m-1\\
t & \mbox{if } k=m.
           \end{cases}\]
Now define $\phi_n:[0,\infty)\to L^1(X)$ by
\[\phi_n(t)=\Phi(\mathbf{s}(n,t),f)\]
when $t>0$ and $\phi_n(0)=f$. By the definition of $S_t$, $|\phi(t)-\phi_n(t)|\to0$ as $n\to\infty$ for each $t\geq0$.

Our task is to demonstrate that $\phi_n$ is measurable for each $n$ in $\N$. Note that $\phi_n(t)$, when $t$ is restricted to the interval $[k2^{-n},(k+1)2^{-n})$, is of the form
\[B_{k,n}\bT_{t-k2^{-n}}f\]
where $B_{k,n}$ is a contraction on $L^1(X)$. It follows that if $E$ is an open set of $L^1(X)$ then
\[\phi_n^{-1}(E)=\bigcup_{k=1}^{\infty}\big\{t\in[k2^{-n},(k+1)2^{-n}):B_{k,n}\bT_{t-k2^{-n}}f\in E\big\}.\]
Hence if the map $\varphi:[0,2^{-n})\to L^1(X)$, defined by
\[\varphi(t)=\bT_tf,\]
is measurable then $\phi^{-1}_n(E)$ can be written as a countable union of measurable sets and consequently $\phi_n$ is measurable. But by Lemma \ref{lem:linear modulus} (d) there is a sequence $\{f_j\}_{j\in\N}$ in $L_+^1(X)$ such that $|\bT_tf-T_tf_j|\to0$ as $j\to\infty$. In other words, $\varphi$ is the pointwise limit of a sequence $\{\varphi_j\}_{j\in\N}$ of measurable functions, defined by $\varphi_j(t)=T_tf_j$, and hence $\varphi$ is measurable. This completes the proof of Theorem \ref{th:contraction implies subpositive}.

\medskip

In preparation for the next section, we present the following lemma. Its proof will be omitted.

\begin{lemma}\label{lem:strong continuity for all p}
Suppose that $\{T_t:t\geq0\}$ is a semigroup on $L^2(X)$ satisfying the contraction property. If $\{T_t:t\geq0\}$ is strongly continuous on $L^p(X)$ for some $p$ in $[1,\infty)$ then $\{T_t:t\geq0\}$ is strongly continuous on $L^q(X)$ for all $q$ in $(1,\infty)$. 
\end{lemma}

%

\section{The Hopf--Dunford--Schwartz ergodic theorem}\label{s:HDS}

We now obtain a vector-valued version of the Hopf--Dunford--Schwartz ergodic theorem for use in Section \ref{s:maximal theorem}. If $\mathcal{T}$ is a bounded strongly measurable semigroup $\{T_s:s\geq0\}$ on $L^p(X)$ then define the operator $A(\mathcal{T},t)$, for positive $t$, by the formula
\[A(\mathcal{T},t)f=\frac{1}{t}\int_0^{t}T_sf\,\dd s\qquad\forall f\in L^p(X).\]
For $f$ in $L^p(X)$, we now define a maximal ergodic function $\mathcal{A}^{\mathcal{T}}f$ by
\begin{equation}\label{eq:HDS maximal function}
\mathcal{A}^{\mathcal{T}}f=\sup_{t>0}|A(\mathcal{T},t)f|.
\end{equation}
A simplified version of the classical Hopf--Dunford--Schwartz ergodic theorem may be stated as follows.

\begin{theorem}\cite[Theorem VIII.7.7]{DS58}\label{th:HDS}
Suppose that $\{T_t:t\geq0\}$ is a measurable semigroup on $L^2(X)$ satisfying the contraction property. Assume that $p\in(1,\infty)$ and denote $\{T_t:t\geq0\}$ by $\mathcal T$. Then the maximal ergodic function operator $\mathcal{A}^{\mathcal{T}}$ satisfies the inequality
\[\norm{\mathcal{A}^{\mathcal T}f}_p\leq2\left(\frac{p}{p-1}\right)^{1/p}\norm{f}_p\qquad\forall f\in L^p(X).\]
\end{theorem}

We will now develop a vector-valued version of this theorem. Fix $p$ in the interval $(1,\infty)$. Suppose that $\mathcal T$ is a strongly continuous semigroup $\{T_t:t\geq0\}$ on $L^2(X)$ satisfying the contraction property. By Lemma \ref{lem:strong continuity for all p}, the semigroup $\mathcal T$ is a strongly continuous semigroup of contractions when viewed as acting on $L^p(X)$. We first show that the bounded linear operator $A(\mathcal{T},t)$ on $L^p(X)$ has an extension to $L^p(X,\B)$ for all positive $t$. By Theorem \ref{th:contraction implies subpositive} there is a measurable semigroup $\{S_t:t\geq0\}$ of positive contractions on $L^p(X)$, which we denote by $\mathcal S$, dominating $\mathcal T$ on $L^p(X)$. Hence $A(\mathcal{S},t)$ is also a positive contraction on $L^p(X)$ for each positive $t$. Moreover,
\begin{equation}\label{eq:A(T,t)f bound}
|A(\mathcal{T},t)f|\leq\frac{1}{t}\int^t_0|T_sf|\,\dd s\leq\frac{1}{t}\int^t_0S_s|f|\,\dd s=A(\mathcal{S},t)|f|
\end{equation}
whenever $f\in L^p(X)$. It follows that $A(\mathcal{T},t)$ has a tensor product extension to $L^p(X,\B)$ for all positive $t$ (see Lemma \ref{lem:subpositivity characterisations}). We can now define a maximal ergodic function operator $\mathcal{A}^{\mathcal T}_{\B}$ by the formula
\begin{equation}\label{eq:HDF vv maximal operator}
\mathcal{A}^{\mathcal T}_{\B}F=\sup_{t>0}|\wt{A}(\mathcal{T},t)F|_{\B}\qquad\forall F\in L^p(X,\B).
\end{equation}

Moreover, if $F\in L^p(X,\B)$ then $\mathcal{A}^{\mathcal T}_{\B}F$ is measurable. To see this, observe that the mapping $t\mapsto A(\mathcal{T},t)f$ is continuous from $(0,\infty)$ to $L^p(X)$ and $\wt{A}(\mathcal{T},t)$ is a contraction on $L^p(X)$ for every $t$ by (\ref{eq:A(T,t)f bound}). Hence the vector-valued mapping $t\mapsto \wt{A}(\mathcal{T},t)F$ is continuous from $(0,\infty)$ to $L^p(X,\B)$ by a simple modification of Lemma \ref{lem:strong continuity}. This implies that the maping $t\mapsto |\wt{A}(\mathcal{T},t)F|_{\B}$ is continuous from $(0,\infty)$ to $L^p(X)$. Therefore the measurable function $\sup_{t\in\Q^+}|\wt{A}(\mathcal{T},t)F|_{\B}$, where $\Q^+$ denotes the set of positive rationals, coincides with $\sup_{t>0}|\wt{A}(\mathcal{T},t)F|_{\B}$.

\begin{corollary}\label{cor:HDS}
Suppose that $\B$ is a Banach space and that $\mathcal{T}$ is a strongly continuous semigroup $\{T_t:t\geq0\}$ on $L^2(X)$ with the contraction property. If $1<p<\infty$ then the maximal ergodic function operator $\mathcal{A}^{\mathcal T}_{\B}$, defined by (\ref{eq:HDF vv maximal operator}), satisfies the inequality
\[\norm{\mathcal{A}_{\B}^{\mathcal T}F}_{L^p(X)}\leq2\left(\frac{p}{p-1}\right)^{1/p}\norm{F}_{L^p(X,\B)}
\qquad\forall f\in L^p(X,\B).\]
\end{corollary}

\begin{proof} Fix $p$ in $(1,\infty)$ and let $\mathcal S$ denote the semigroup dominating $\mathcal T$ that was introduced above.  For $F$ in $L^p(X,\B)$,
\begin{align*}
\mathcal{A}^{\mathcal{T}}_{\B}F
&=\sup_{t>0}|\wt{A}(\mathcal{T},t)F|_{\B}\\
&\leq \sup_{t>0}A(\mathcal{S},t)|F|_{\B}\\
&=\mathcal{A}^{\mathcal S}|F|_{\B}.
\end{align*}
The result follows upon taking the $L^p(X)$ norm of both sides and applying Theorem \ref{th:HDS}.
\end{proof}

\section{The vector-valued maximal theorem}\label{s:maximal theorem}

The main result of this section is a vector-valued version of Theorem \ref{th:scalar maximal}. It gives an $L^p$ estimate for the maximum function $\mathcal{M}^{\psi}_{\B}F$ (defined by (\ref{eq:def of maximal function})) under the assumption that the generator $-\wt{L}$ of $\{\wt{T}_t:t\geq0\}$ has bounded imaginary powers. 
 
\begin{theorem}\label{th:maximal} Suppose that $(X,\mu)$ is a $\sigma$-finite measure space, $\B$ is a Banach space, $1<p<\infty$ and $\{T_t:t\geq0\}$ is a strongly continuous semigroup on $L^2(X)$ with the contraction property. If there exists $\omega$ less than $\pi/2-\psi$ and a positive constant $K$ such that $\wt{L}$ has bounded imaginary powers satisfying the norm estimate
\begin{equation}\label{eq:BIP}
\Vert\wt{L}^{iu}F\Vert_{L^p(X,\B)} \leq Ke^{\omega|u|}\norm{F}_{L^p(X,\B)}
\qquad\forall F\in L^p(X,\B)\quad\forall u\in\R,
\end{equation}
then $\{\wt{T}_t:t\geq0\}$ has a bounded analytic continuation in $L^p(X,\B)$ to the sector $\Gamma_{\psi}$ and there is a constant $C$ such that the maximal function operator $\mathcal{M}^{\psi}_{\B}$ satisfies the inequality
\begin{equation}\label{eq:maximal function bound}
\norm{\mathcal{M}^{\psi}_{\B}F}_{L^p(X)}\leq C\norm{F}_{L^p(X,\B)}\qquad\forall F\in L^p(X,\B).
\end{equation}
\end{theorem}

\begin{proof} Assume the hypotheses of the theorem. Since $\wt{L}$ has bounded imaginary powers satisfying (\ref{eq:BIP}), $-\wt{L}$ generates a uniformly bounded semigroup on $L^p(X,\B)$ with analytic continuation to any sector $\Gamma_{\psi_0}$, where
\[\psi_0<\frac{\pi}{2}-\omega,\]
by a result of J. Pr\"uss and H. Sohr \cite[Theorem 2]{PS90}. Hence the operator $\mathcal{M}^{\psi}_{\B}$ is well-defined. It remains to show (\ref{eq:maximal function bound}).

Take $F$ in $L^p(X,\B)$ and $z$ in $\overline{\Gamma}_{\psi}\bs\{0\}$. Write $z$ as $e^{i\theta}t$, where $|\theta|\leq\psi$ and $t>0$. The key idea of the proof is to decompose $\wt{T}_zF$ into two parts:
\begin{equation}\label{eq:decomposition3}
\wt{T}_zF=\frac{1}{t}\int_0^te^{-s\wt{L}}F\,\dd s+\Big[e^{-z\wt{L}}F-\frac{1}{t}\int_0^te^{-s\wt{L}}F\,\dd s\Big].
\end{equation}
Define the function $m_{\theta}$ on $(0,\infty)$ by
\begin{equation}\label{eq:m_theta}
m_{\theta}(\lambda)=\exp(-e^{i\theta}\lambda)-\int_0^1e^{-s\lambda}\,\dd s\qquad\forall\lambda>0.
\end{equation}
Then (\ref{eq:decomposition3}) can be rewritten formally as
\[\wt{T}_zF=\frac{1}{t}\int_0^te^{-s\wt{L}}F\,\dd s+m_{\theta}(t\wt{L})F,\]
whence
\[\sup_{z\in\overline{\Gamma}_{\psi}\bs\{0\}}|\wt{T}_zF|_{\B}
\leq\sup_{t>0}\left|\frac{1}{t}\int_0^te^{-s\wt{L}}F\,\dd s\right|_{\B}
+\sup_{t>0}\;\sup_{|\theta|\leq\psi}|m_{\theta}(t\wt{L})F|_{\B}.\]
If we take the $L^p(X)$ norm of both sides then we have, formally at least,
\begin{equation}\label{eq:decompsition4}
\norm{\mathcal{M}^{\psi}_{\B}F}_p
\leq\norm{\mathcal{A}_{\B}^{\mathcal T}F}_p+\norm{\sup_{t>0}\;\sup_{|\theta|\leq\psi}|m_{\theta}(t\wt{L})F|_{\B}}_p,
\end{equation}
where $\mathcal T$ denotes the semigroup $\{T_t:t\geq0\}$ and $\mathcal{A}^{\mathcal T}_{\B}$ is the operator defined by (\ref{eq:HDF vv maximal operator}).
By Corollary \ref{cor:HDS}, the first term on the right-hand side is majorised by $2[p/(1-p)]^{1/p}\norm{F}_{L^p(X,\B)}$. We need to control the second term.

Write $n_{\theta}$ for $m_{\theta}\circ \exp$ and observe that
\begin{equation}\label{eq:m and n}
m_{\theta}(\lambda)=\frac{1}{2\pi}\int_{-\infty}^{\infty}\hat{n}_{\theta}(u)\lambda^{iu}\,\dd u,
\end{equation}
where $\hat{n}_{\theta}$ denotes the Fourier transform of $n_{\theta}$.
Calculation using complex analysis shows that
\[\hat{n}_{\theta}(u)=\big(e^{-\theta u}-(1+iu)^{-1}\big)\Gamma(iu)\qquad\forall u\in\R,\]
and the theory of the $\Gamma$-function (see, for example, \cite[p.~151]{eT78}) gives the estimate
\[|\hat{n}_{\theta}(u)|\leq C_0\exp\big((|\theta|-\pi/2)|u|\big)\qquad\forall u\in\R,\]
where $C_0$ is a constant independent of $u$ and $\theta$. Thus, the existence of bounded imaginary powers of $\wt{L}$ gives
\begin{align*}
\sup_{t>0}\,\sup_{|\theta|\leq\psi}|m_{\theta}(t\wt{L})F|_{\B}
&\leq\sup_{t>0}\,\sup_{|\theta|\leq\psi}
\frac{1}{2\pi}\int_{-\infty}^{\infty}|\hat{n}_{\theta}(u)|\,|(t\wt{L})^{iu}F|_{\B}\,\dd u\\
&\leq \sup_{t>0}\,\sup_{|\theta|\leq\psi}
\frac{1}{2\pi}\int_{-\infty}^{\infty}C_0e^{(|\theta|-\pi/2)|u|}|t^{iu}|\,|\wt{L}^{iu}F|_{\B}\,\dd u\\
&\leq\frac{C_0}{2\pi}\int_{-\infty}^{\infty}e^{(\psi-\pi/2)|u|}|\wt{L}^{iu}F|_{\B}\,\dd u.
\end{align*}
Taking the $L^p(X)$ norm of both sides of the above inequality and applying (\ref{eq:BIP}) gives
\begin{align*}
\norm{\sup_{t>0}\;\sup_{|\theta|\leq\psi}\;|m_{\theta}(t\wt{L})F|_{\B}}_p
&\leq\frac{C_0}{2\pi}\int_{-\infty}^{\infty}e^{(\psi-\pi/2)|u|}\norm{\wt{L}^{iu}F}_p\,\dd u\\
&\leq \frac{C_0K}{2\pi}\int_{-\infty}^{\infty}e^{(\psi-\pi/2)|u|}
e^{\omega|u|}\norm{F}_p\dd u\\
&<C_1\norm{F}_{L^p(X,\B)},
\end{align*}
since $\psi-\pi/2+\omega<0$, and where $C_1$ is a positive constant independent of $F$. Now (\ref{eq:decompsition4}) and Corollary \ref{cor:HDS} yields (\ref{eq:maximal function bound}) for some positive constant $C$.

The opening formal calculations can be justified by working backwards, provided that the function
\begin{equation}\label{fun:m theta}
\sup_{t>0}\;\sup_{|\theta|\leq\psi}|m_{\theta}(t\wt{L})F|_{\B}
\end{equation}
is measurable. Since the map $z\mapsto\wt{T}_zF$ is continuous from $\overline{\Gamma}_{\psi}$ to $L^p(X,\B)$, the map $(t,\theta)\mapsto|m_{\theta}(t\wt{L})F|_{\B}$ is continuous from $(0,\infty)\times[-\psi,\psi]$ to $L^p(X)$. Hence
\[\sup_{t>0}\;\sup_{|\theta|\leq\psi}|m_{\theta}(t\wt{L})F|_{\B}
=\sup_{(t,\theta)\in R}|m_{\theta}(t\wt{L})F|_{\B},\]
where $R$ is the denumerable set $\big((0,\infty)\times[-\psi,\psi]\big)\cap\Q^2$. Since each $m_{\theta}(t\wt{L})F$ is measurable in $L^p(X,\B)$ it follows that (\ref{fun:m theta}) is measurable in $L^p(X)$. This completes the proof.
\end{proof}

\section{Bounded imaginary powers of the generator}\label{s:BIP}

In this section we examine circumstances under which the bounded imaginary power estimate (\ref{eq:BIP}), one of the hypotheses of the preceding theorem and corollary, is satisfied. A fruitful (and in our context, necessary) setting  is when the Banach space $\B$ has the \textit{UMD property}. A Banach space $\B$ is said to be a \textit{UMD space} if one of the following equivalent statements hold:
\begin{enumerate}
\item[(a)] The Hilbert transform is bounded on $L^p(X,\B)$ for one (and hence all) $p$ in $(1,\infty)$.
\item[(b)] If $1<p<\infty$ then $\B$-valued martingale difference sequences on $L^p(X,\B)$ converge unconditionally.
\item[(c)] If $1<p<\infty$ then $(-\Delta)^{iu}\otimes I_{\B}$ extends to a bounded operator on $L^p(\R,\B)$ for every $u$ in $\R$ (a result due to S. Guerre-Delabri{\`e}re \cite{GD91}).
\end{enumerate}
Several other characterisations of UMD spaces exist (see, for example, \cite{dB81} and the survey in \cite{RdF86}) but those cited here are, for different reasons, the most relevant to our discussion. If the Hilbert transform, which corresponds to the multiplier function $u\mapsto i\,\mathrm{sgn}(u)$, is bounded on $L^p(X,\B)$ then one can establish vector-valued versions of some Fourier multiplier theorems (such as Mikhlin's multiplier theorem \cite{fZ89}). This fact is used below to establish Theorem \ref{th:H infty functional calculus}. The second characterisation gave rise to the name UMD. The third characterisation shows that, in general, $\B$ must be a UMD space if $\wt{L}$ is to have bounded imaginary powers, since $-\Delta$ generates the Gaussian semigroup. Examples of UMD spaces include, when $1<p<\infty$, the classical $L^p(X)$ spaces and the Schatten--von Neumann ideals $\mathcal{C}^p$. Moreover, if $\B$ is a UMD space then its dual $\B^*$, closed subspaces of $\B$, quotient spaces of $\B$ and $L^p(X,\B)$ when $1<p<\infty$ also inherit the UMD property.

It was shown by Hieber and Pr\"uss \cite{HP98} that when $1<q<\infty$ the generator of a UMD-valued extension of a bounded strongly continuous positive semigroup on $L^q(X)$ has a bounded $H^{\infty}$-functional calculus. The next result says that the same is true if the positivity condition is relaxed to subpositivity (assuming that the UMD-valued extension of the semigroup is bounded), though it is convenient in the present context to state it for semigroups possessing the contraction property. First we introduce some notation. If $\sigma\in(0,\pi]$ then let $H^{\infty}(\Gamma_{\sigma})$ denote the Banach space of all bounded analytic functions defined on $\Gamma_{\sigma}$ with norm \[\norm{f}_{H^{\infty}(\Gamma_{\sigma})}=\sup_{z\in\Gamma_{\sigma}}|f(z)|.\]

\begin{theorem}\label{th:H infty functional calculus}
Suppose that $1<q<\infty$ and $\B$ is a UMD space. If $\{T_t:t\geq0\}$ is a strongly continuous semigroup on $L^2(X)$ satisfying the contraction property and $-\wt{L}$ is the generator of its tensor extension $\{\wt{T}_t:t\geq0\}$ to $L^q(X,\B)$, then $\wt{L}$ has a bounded $H^{\infty}(\Gamma_{\sigma})$-calculus for all $\sigma$ in $(\pi/2,\pi]$. Consequently, for every $\sigma\in(\pi/2,\pi]$ there exists a positive constant $C_{q,\sigma}$ such that
\begin{equation}\label{eq:BIP rough estimate}
\Vert\wt{L}^{iu}F\Vert_{L^q(X,\B)}\leq C_{q,\sigma}e^{\sigma|u|}\norm{F}_{L^q(X,\B)}\qquad
\forall F\in L^q(X,\B)\quad\forall u\in\R.
\end{equation}
\end{theorem}

\begin{proof} Since the semigroup $\{T_t:t\geq0\}$ can be extended to a subpositive strongly continuous semigroup of contractions on $L^q(X)$, it has a dilation to a bounded $c_0$-group on $L^q(X')$ for some measure space $(X',\mu')$. In other words, there exists a measure space $(X',\mu')$, a strongly continuous group $\{U_t:t\in\R\}$ of subpositive contractions on $L^q(X')$, a positive isometric embedding $D:L^q(X)\to L^q(X')$ and a subpositive contractive projection $P:L^q(X')\to L^q(X')$ such that
\[DT_t=PU_tD\qquad\forall t\geq0\]
(see the result of G. Fendler \cite[pp.~737--738]{gF97} which extends the work of Coifman, Rochberg, and Weiss \cite{CRW77}). Lifting this identity to its $\B$-valued extension, we see that the semigroup $\{\wt{T}_t:t\geq0\}$ on $L^q(X,\B)$ has a dilation to a bounded $c_0$-group $\{\wt{U}_t:t\in\R\}$ on $L^q(X',\B)$.

Let $-\wt{L}$ denote the generator of $\{\wt{T}_t:t\geq0\}$. Then the dilation implies that $\wt{L}$ has a bounded $H^{\infty}(\Gamma_{\sigma})$-calculus for all $\sigma$ in $(\pi/2,\pi]$ (see \cite{HP98} or the exposition in \cite[pp.~212--214]{KW04}, where the $H^{\infty}$-calculus is first constructed for the generator of the group $\{\wt{U}_t:t\in\R\}$ using the vector-valued Mikhlin multiplier theorem in conjunction with the transference principle, and then projected back to the generator $-\wt{L}$ of $\{\wt{T}_t:t\geq0\}$ via the dilation).

The bounded $H^{\infty}(\Gamma_{\sigma})$-calculus gives a positive constant $C_{q,\sigma}$ such that
\[\Vert f(\wt{L})\Vert_{L^q(X,\B)}\leq C_{q,\sigma}\norm{f}_{H^{\infty}(\Gamma_{\sigma})}
\qquad\forall f\in H^{\infty}(\Gamma_{\sigma}).\]
If $f(z)=z^{iu}$ for $u$ in $\R$ then (\ref{eq:BIP rough estimate}) follows.
\end{proof}

The theorem above suggests that the problem of finding bounded imaginary powers of $\wt{L}$ is critical to $L^2(X,\B)$. That is, if
\[\Vert\wt{L}^{iu}F\Vert_{L^2(X,\B)}\leq Ce^{\omega|u|}\norm{F}_{L^2(X,\B)} \qquad\forall F\in L^2(X,\B)\:\forall u\in\R\]
for some $\omega$ less than $\pi/2-\psi$ then one could interpolate between the $L^2$ estimate and (\ref{eq:BIP rough estimate}) to obtain (\ref{eq:BIP}). Unfortunately, suitable $L^2(X,\B)$ bounded imaginary power estimates, where $\B$ is a nontrivial UMD space, appear to be absent in the literature, even when $\wt{L}$ is the Laplacian. However, if $\B$ is a Hilbert space, such estimates are available via spectral theory.

\begin{lemma} Suppose that $\mathcal{H}$ is a Hilbert space. If $\{T_t:t\geq0\}$ is a symmetric diffusion semigroup on $L^2(X)$ then the generator $-\wt{L}$ of the $\Hi$-valued extension $\{\wt{T}_t:t\geq0\}$ to $L^2(X,\Hi)$ satisfies
\begin{equation}\label{eq:BIP L^2(X,H)}
\Vert\wt{L}^{iu}F\Vert_{L^2(X,\Hi)}\leq \norm{F}_{L^2(X,\Hi)}\qquad\forall F\in L^2(X,\Hi)\quad\forall u\in\R.
\end{equation}
\end{lemma}

\begin{proof} It is not hard to check that the tensor product extension to $L^2(X,\Hi)$ of the semigroup $\{T_t:t\geq0\}$ is a semigroup of selfadjoint contractions on $L^2(X,\Hi)$. Its generator $-\wt{L}$ is therefore selfadjoint on $L^2(X,\Hi)$  and hence $\wt{L}$ has nonnegative spectrum. Spectral theory now gives estimate (\ref{eq:BIP L^2(X,H)}).
\end{proof}

To obtain (\ref{eq:BIP}) we shall interpolate between (\ref{eq:BIP rough estimate}) and (\ref{eq:BIP L^2(X,H)}). Hence we consider the class of UMD spaces whose members $\B$ are isomorphic to closed subquotients of a complex interpolation space $(\Hi,\U)_{[\theta]}$, where $\Hi$ is a Hilbert space, $\U$ is a UMD space and $0<\theta<1$. Members of this class include the UMD function lattices on a $\sigma$-finite measure space (such as the reflexive $L^p(X)$ spaces) by a result of Rubio de Francia (see \cite[Corollary, p.~216]{RdF86}), the reflexive Sobolev spaces (which are subspaces of products of $L^p$ spaces) and the reflexive Schatten--von Neumann ideals. This class can be further extended to include many operator ideals by combining Rubio de Francia's theorem with results due to P. Dodds, T. Dodds and B. de Pagter \cite{DDdP92} which show that the interpolation properties of noncommutative spaces coincide with those of their commutative counterparts under fairly general conditions. It was asked in \cite{RdF86} whether the described class of UMD spaces includes all UMD spaces. It appears that this is still an open question.

\begin{corollary}\label{cor:BIP} Suppose that $\B$ is a UMD space isomorphic to a closed subquotient of a complex interpolation space $(\Hi,\U)_{[\theta]}$, where $\Hi$ is a Hilbert space, $\U$ is a UMD space and $0<\theta<1$. Suppose also that $\{T_t:t\geq0\}$ is a symmetric diffusion semigroup on $L^2(X)$ and denote by $-\wt{L}$ the generator of its tensor extension to $L^p(X,\B)$, where $1<p<\infty$. If
\begin{equation}\label{eq:p and theta}
|2/p-1|<\theta
\end{equation}
and
\[0\leq\psi<\frac{\pi}{2}(1-\theta)\]
then there exists $\omega$ less than $\pi/2-\psi$ such that $\wt{L}$ has bounded imaginary powers on $L^p(X,\B)$ satisfying estimate (\ref{eq:BIP}).
\end{corollary}

\begin{proof}
Assume the hypotheses of the corollary. Note that
\[\frac{\pi}{2}<\frac{1}{\theta}\left(\frac{\pi}{2}-\psi\right)\]
so that if $\sigma$ is the arithmetic mean of $\pi/2$ and $(\pi/2-\psi)/\theta$ then $\sigma>\pi/2$ and $\sigma\theta<\pi/2-\psi$. Now choose $q$ such that
\[\frac{1}{p}=\frac{1-\theta}{2}+\frac{\theta}{q}.\]
Inequality (\ref{eq:p and theta}) guarantees that $1<q<\infty$. Interpolating between (\ref{eq:BIP L^2(X,H)}) and (\ref{eq:BIP rough estimate}) (for the space $L^q(X,\U)$) gives
\[\Vert\wt{L}^{iu}\Vert_{L^p(X,\B)}\leq C_{q,\sigma}^{\theta}e^{\sigma\theta|u|}\norm{F}_{L^p(X,\B)}
\qquad\forall F\in L^p(X,\B)\quad\forall u\in\R.\]
If $\omega=\sigma\theta$ then (\ref{eq:BIP}) follows, completing the proof.
\end{proof}

\section{Proof of Theorem \ref{th:main semigroup theorem}}\label{s:proof of main result}

In this final section we complete the proof of Theorem \ref{th:main semigroup theorem}. Suppose the hypotheses of the theorem. Parts (a) and (b) follow immediately from Theorem \ref{th:maximal} and Corollary \ref{cor:BIP}. Part (c) will be deduced from the vector-valued maximal theorem and the pointwise convergence of $\{T_t:t\geq0\}$ (see Corollary \ref{cor:sv ptwise}).

For ease of notation, write $z\to0$ as shorthand for $z\to0$ with $z$ in $\overline{\Gamma}_{\psi}$.

Suppose that $F\in L^p(X,\B)$ and $\epsilon>0$. There exists a function $G$ in $L^p(X)\otimes\B$ such that $\norm{G-F}_{L^p(X,\B)}<\epsilon$. Write $G$ as $\sum_{k=1}^nu_kf_k$, where $n$ is a positive integer, $\{u_k\}_{k=1}^n$ is contained in $\B$ and $\{f_k\}_{k=1}^n$ is contained in $L^p(X)$.
Hence, for almost every $x$ in $X$,
\begin{align*}
\limsup_{z\to0}|\wt{T}_zF(x)-F(x)|_{\B}
&\leq \limsup_{z\to0}|\wt{T}_zF(x)-\wt{T}_zG(x)|_{\B}+|G(x)-F(x)|_{\B}\\
&\quad+\limsup_{z\to0}|\wt{T}_zG(x)-G(x)|_{\B}\\
&\leq \sup_{z\in\overline{\Gamma}_{\psi}}|\wt{T}_z(F-G)(x)|_{\B}+|G(x)-F(x)|_{\B}\\
&\quad+\sum_{k=1}^n\big|u_k\big|_{\B}\limsup_{z\to0}\big|T_zf_k(x)-f_k(x)\big|\\
&\leq 2\mathcal{M}^{\psi}_{\B}(G-F)(x),
\end{align*}
since Corollary \ref{cor:sv ptwise} implies that
\[\lim_{z\to0}|T_zf_k(x)-f_k(x)|=0\]
for each $k$ and for almost every $x$ in $X$. By taking the $L^p(X)$ norm and applying Theorem \ref{th:maximal} we obtain
\[\norm{\limsup_{z\to0}|\wt{T}_zF-F|_{\B}}_p\leq 2\norm{\mathcal{M}^{\psi}_{\B}(G-F)}_p<2C\epsilon,\]
where the positive constant $C$ is independent of $F$ and $G$.
Since $\epsilon$ is an arbitrary positive number,
\[\limsup_{z\to0}|\wt{T}_zF(x)-F(x)|_{\B}=0\]
for almost every $x$ in $X$, proving the theorem.

\vspace*{.3cm}

\noindent\textit{Acknowledgements.}\quad
It is a pleasure to acknowledge the contribution of Michael Cowling, who introduced me to this problem and gave many helpful comments and suggestions. Thanks also go to Ian Doust, who helped shed light on a couple of technical obstacles, and to Pierre Portal and Tuomas Hyt\"onen for the interest they showed in this work. The referee made several valuable comments, one which strengthened the main result of this paper. This research was funded by an Australian Postgraduate Award and by the Australian Research Council's Centre of Excellence for Mathematics and Statistics of Complex Systems.

%
%

\end{document}